\documentclass[12pt]{amsart}
\usepackage{amssymb}
\usepackage{amsmath}
\usepackage{amsthm}
\newcommand{\N}{{\mathcal N}}
\newcommand{\C}{\ensuremath{\mathbb{C}}}
\newcommand{\R}{\ensuremath{\mathbb{R}}}
\newcommand{\vk}{\varkappa}
\newcommand{\p}{\partial}
\newcommand{\s}{{\rm Symb}}
\newcommand{\T}{{T^\ast M}}
\newcommand{\D}{{\mathcal D}}
\newcommand{\W}{{\mathcal W}}
\newtheorem{lemma}{Lemma}
\newtheorem{proposition}{Proposition}
\newtheorem{theorem}{Theorem}
\newtheorem{corollary}{Corollary}
\begin{document}

\title{On Dequantization of Fedosov's Deformation Quantization}

\author[A.V. Karabegov]{Alexander V. Karabegov}\footnote{Research was partially supported by an ACU Math/Science grant.}
\address[Alexander V. Karabegov]{Department of Mathematics and Computer Science, Abilene Christian University, ACU Box 28012, 252 Foster Science Building, Abilene, TX 79699-8012}
\email{alexander.karabegov@math.acu.edu}

\begin{abstract}
 To each natural deformation quantization on a Poisson manifold $M$ we associate a Poisson morphism from the formal neighborhood of the zero section of $\T$ to the formal neighborhood of the diagonal of the product $M\times \widetilde M$, where $\widetilde M$ is a copy of $M$ with the opposite Poisson structure. We call it dequantization of the natural deformation quantization. Then we ``dequantize" Fedosov's quantization.
\end{abstract}
\subjclass[2000]{Primary: 53D55; Secondary:  53D05}
\keywords{deformation quantization, dequantization, Poisson morphism}
\date{Received 30 July 2003}
\maketitle

\section{Introduction}

Quantization relates to geometric structures of classical mechanics their algebraic quantum counterparts.  However,  between the classical and quantum worlds there is a world of  ``semiclassical" objects that are still of geometric nature but  which inherit certain properties of quantum structures. There is a large and important class of such semiclassical objects, symplectic groupoids, introduced independently by Karas\"ev  \cite{Ka},  Weinstein \cite{W}, and Zakrzewski \cite{Z}. Symplectic groupoids were extensively studied as geometric objects (see \cite{CrF} and references therein). In \cite{CF2} Cattaneo and Felder constructed a symplectic groupoid based on the Poisson sigma model whose perturbative quantization yields the Kontsevich star-product (see \cite{CF1}). In what follows we associate to each natural  star-product (see the definition in Section \ref{S:nat}) on a Poisson manifold $M$ a formal geometric object, $(\T,Z)$, the formal neighborhood of the zero section $Z$ of the cotangent bundle $\T$ endowed with the ``source" and ``target"  mappings to $M$ which are a Poisson and anti-Poisson morphisms  respectively (we give the necessary terminology related to the formal neighborhoods in the Appendix).  This object can be thought of as a ``formal symplectic groupoid" at least if $M$ is symplectic.
We call it ``dequantization" of the initial star-product. The goal of this paper is to calculate the source and target mappings associated to a Fedosov star-product. These mappings can be described independently in terms of the (formal)  geometric structure of Fedosov's quantization revealed by Emmrich and Weinstein in \cite{EW}.

{\bf Acknowledgments.}  We are very grateful to Hovhannes Khudaverdian, Simon Lyakhovich, and Boris Tsygan for stimulating discussions and to Abilene Christian University for supporting our research.
\bigskip

\section{Natural Deformation Quantizations}\label{S:nat}

\bigskip

Let $(M,\{\cdot,\cdot\})$ be a Poisson manifold. Denote
by $C^\infty(M)[[\nu]]$ the space of formal series in $\nu$ with
coefficients from $C^\infty(M)$. As introduced in \cite{BFFLS}, a formal
differentiable deformation quantization on $M$ is an associative algebra
structure on $C^\infty(M)[[\nu]]$ with the $\nu$-linear and $\nu$-adically
continuous product $\ast$ (named star-product) given on $f,g\in
C^\infty(M)$ by the formula
\begin{equation} \label{E:star} 
f \ast g = \sum_{r = 0}^\infty \nu^r C_r(f,g), 
\end{equation} 
where $C_r,\, r\geq 0,$ are bidifferential operators on $M$, $C_0(f,g) =
fg$ and $C_1(f,g) - C_1(g,f) = i\{f,g\}$.  We adopt the convention that
the unit of a star-product is the unit constant. Two differentiable
star-products $\ast,\ast'$ on a Poisson manifold $(M,\{\cdot,\cdot\})$ are
called equivalent if there exists an isomorphism of algebras $B:
(C^\infty(M)[[\nu]],\ast) \to (C^\infty(M)[[\nu]],\ast')$ of the form $B=
1 + \nu B_1 + \nu^2 B_2 + \dots,$ where $B_r, r\geq 1,$ are differential
operators on $M$. The existence and classification problem for deformation quantization was first solved in the non-degenerate (symplectic) case (see \cite{DWL}, \cite{OMY}, \cite{F1} for existence proofs and \cite{F2}, \cite{NT}, \cite{D}, \cite{BCG}, \cite{X} for classification) and then Kontsevich \cite{K} showed that every Poisson manifold admits a deformation quantization and that the equivalence classes of deformation quantizations can be parameterized by the formal deformations of the Poisson structure. 

It turns out that all the explicit constructions of star-products enjoy the following property: for all $r \geq 0$ the bidifferential operator $C_r$ in (\ref{E:star}) is of order not greater than $r$ in both arguments (the most important examples are Fedosov star-products on symplectic manifolds and  the Kontsevich star-product on $\R^n$ endowed with an arbitrary Poisson bracket). The star-products with this property were called natural by Gutt and Rawnsley in \cite{GR}, where general properties of such star-products were studied (before these star-products were said to be of Vey type).

{\it Remark.} In this Letter the word "formal" is used very often in different contexts. In order to avoid unnecessary confusion, from now on we will call the objects that are formal series in the formal parameter $\nu$ $\nu$-formal.  
\section{Dequantization of a natural deformation quantization}

Let $\D = \D(M)$ denote the algebra of differential operators with smooth complex-valued coefficients on a real manifold $M$. It has a natural filtration with respect to the order of differential operators, $0 \subset D_1 \subset D_2 \subset \dots$. The symbol mapping $\s_r: \D_r \to C^\infty(T^\ast M)$ maps the operators of order not greater than $r$ to the functions on $\T$ that are homogeneous polynomials of degree $r$ on the fibres. The kernel of $\s_r$ is $\D_{r-1}$. If $\{x^k, \xi_l\}$ are local  coordinates on $\T$ (as usual, $\{\xi_l\}$ are fibre coordinates dual to the base coordinates $\{x^k\}$) and an operator $X\in\D_r$ has the leading term $X^J \p_J$ where $J$ runs over the multi-indices of order $r$ and $\p_k = \p/ \p x^k$, then $\s_r (X) = X^J \xi_J$. Denote by $\{\cdot,\cdot\}_\T$ the natural Poisson bracket on $\T$. In local coordinates 
\[
\{f,g\}_\T = \frac{\p f}{\p \xi_k}\frac{\p g}{\p x^k} -  \frac{\p g}{\p \xi_k}\frac{\p f}{\p x^k}.
\]
 For two arbitrary differential operators $X\in\D_k$ and $Y\in \D_l$ it is well known that $XY \in \D_{k+l},\ [X,Y] \in \D_{k+l-1}$ and
\begin{equation}\label{E:product}
\s_{k+l} XY = \s_k X\cdot \s_l Y,
\end{equation}
\begin{equation}\label{E:comm} 
\s_{k+l-1}([X,Y]) = \{\s_k(X),\s_l(Y)\}_\T.    
\end{equation}

For each point $x\in M$ denote by $\D^x$ the space of all distributions supported at $x$. It is known that $\D^x$ consists of the linear combinations of partial derivatives of the Dirac distribution $\delta_x$. The space $\D^x$ has also a natural filtration $0 \subset \D^x_0 \subset \D^x_1 \subset\dots$ with respect to the order of distributions. The symbol mapping $\s_r: \D^x_r \to \C[\xi]$ maps the distributions of order not greater than $r$ to homogeneous polynomials on $T^\ast_xM$ of degree $r$. Consider the mapping $A \mapsto A_x$ which maps an operator $A$ to the distribution $A_x$ such that $(A_x,\phi) = (A\phi)(x)$, where $\phi$ is a test function. This mapping maps $\D_r$ onto $\D^x_r$.
Moreover, the symbol of $A_x$ is the restriction of the symbol of $A$ to the fibre $T^\ast_xM,\ \s_r(A_x)(\xi) = \s_r(A)(x,\xi)$. Thus the symbol of $A$ is completely determined by the symbols of all $A_x$.

We call a $\nu$-formal differential operator  $A = A_0 + i\nu A_1 + (i\nu)^2 A_2 +\dots$ {\it natural} if the order of $A_r$ is not greater than $r$. 
The natural $\nu$-formal differential operators form an algebra which we denote by $\N$. We define a $\sigma$-symbol of a natural $\nu$-formal differential operator $A = A_0 + i\nu A_1 + (i\nu)^2 A_2 +\dots$ as the formal series $\sigma(A) = \s_0(A_0) + \s_1(A_1) + \dots$. The $\sigma$-symbols are elements of $C^\infty(\T,Z)$ (see the Appendix). Similarly, we call a $\nu$-formal distribution $u = u_0 + i\nu u_1 + \dots$ natural if the order of the distribution $u_r$ is not greater than $r$. A $\sigma$-symbol of $u$ is $\sigma(u) = \s_0(u_0) + \s_1(u_1) +\dots$.

It follows from (\ref{E:product}) that for $A,B \in \N$
\begin{equation}\label{E:sprod}
\sigma (AB) = \sigma(A)\cdot \sigma(B).
\end{equation}
Thus the mapping $\sigma: A \mapsto \sigma(A)$ is a homomorphism from the algebra $\N$ to the algebra $C^\infty(\T,Z)$. Notice that the kernel of the mapping $\sigma$ is $\nu\N$.

It follows from (\ref{E:comm}) that for two natural $\nu$-formal differential operators $A,B \in \N$ the $\nu$-formal differential operator $(1/i\nu)[A,B]$ is natural and 
\begin{equation}\label{E:scomm}
\sigma\left(\frac{1}{i\nu}[A,B]\right) = \{\sigma(A),\sigma(B)\}_\T.     
\end{equation}

Consider a star-product  $\ast$ on a Poisson manifold $(M, \{\cdot,\cdot\})$. For $f,g \in C^\infty(M)[[\nu]]$ denote by $L_f$ and $R_g$ the operators of left $\ast$-multi\-pli\-cation by $f$ and of right $\ast$-multi\-pli\-cation by $g$ respectively, so that $L_f g = f \ast g = R_g$. The associativity of $\ast$ is equivalent to the fact that $[L_f, R_g]=0$.

A star-product $\ast$ on $M$ is natural iff for any  $f,g \in C^\infty(M)[[\nu]]$ the operators $L_f, R_g$ are natural. For a given natural star-product $\ast$ we introduce the following mappings 
$S, T : C^\infty(M)[[\nu]] \to C^\infty(\T,Z)$. For $f \in C^\infty(M)[[\nu]]$ we set $Sf = \sigma(L_f)$ and $Tf = \sigma(R_f)$. Since the kernel of $\sigma$ is $\nu\N$, then for $f = f_0 + \nu f_1 + \dots$ we have that $Sf = Sf_0$ and $Tf=Tf_0$. Now we will show that the mapping $S: C^\infty(M) \to C^\infty(\T,Z)$ is a Poisson morphism. Take $f,g \in C^\infty(M)$.  Since $f\ast g = fg\pmod {\nu}$, we get from (\ref{E:sprod}):
\[
 S(fg) = S(f\ast g) =\sigma(L_{f\ast g}) = \sigma(L_fL_g)=\sigma (L_f)\sigma (L_g) = Sf \cdot Sg.
\]
Since the $\ast$-commutator $[f,g]_\ast = i\nu\{f,g\}\pmod {\nu^2}$, it follows from (\ref{E:scomm}) that 
\begin{eqnarray*}
S(\{f,g\}) = \sigma (L_{\{f,g\}}) = \sigma (L_{(1/i\nu)[f,g]_\ast\}}) = \sigma\left(\frac{1}{i\nu}[L_f,L_g]\right)\\
= \{\sigma(L_f),  \sigma(L_g)\}_\T = \{Sf,Sg\}_\T.
\end{eqnarray*}
Similarly, the mapping $T: C^\infty(M) \to C^\infty(\T,Z)$ is an anti-Poisson morphism, i.e. $T(fg) = Tf \cdot Tg$ and $T(\{f,g\}) = - \{Tf,Tg\}_\T$. Otherwise speaking, $T$ is a Poisson morphism from $C^\infty(\widetilde M)$ to $C^\infty(\T,Z)$, where $\widetilde M$ is a copy of the manifold $M$ endowed with the opposite Poisson bracket $-\{\cdot,\cdot\}$.
Since for any $f,g\in C^\infty(M)$ the operators $L_f$ and $R_g$ commute, we have from (\ref{E:scomm}) that
\[
\{Sf, Tg\}_\T =\{\sigma(L_f),\sigma(R_g)\}_\T = \sigma \left(\frac{1}{i\nu}[L_f,R_g]\right)=0.
\]
Thus we can combine the mappings $S$ and $T$ into a single Poisson morphism $S\otimes T: C^\infty(M\times\widetilde M) \to C^\infty(\T,Z)$ such that $S \otimes T(f \otimes g)= Sf\cdot Tg$. Here $(f\otimes g)(x,\tilde x) = f(x)\cdot g(\tilde x)$.

The operator $C_1$ of a natural star-product $\ast$ can be written for $f,g \in C^\infty(M)$  as 
\begin{equation}\label{E:c1}
C_1(f,g) = \frac{i}{2}\Lambda^{jk}\frac{\p f}{\p x^j}\frac{\p g}{\p x^k},
\end{equation}
where $\Lambda^{jk}$ is a tensor on $M$ such that $\omega^{jk} = \frac{1}{2}\left(\Lambda^{jk} - \Lambda^{kj}\right)$ is the Poisson tensor on $M$. It follows from (\ref{E:c1}) that for $f\in C^\infty(M)$ we have $L_f = f + \frac{i\nu}{2}\Lambda^{jk}\frac{\p f}{\p x^j}\frac{\p }{\p x^k} + \dots$. Thus
$Sf = \sigma(L_f) = f + \frac{1}{2}\Lambda^{jk}\frac{\p f}{\p x^j}\xi_k \pmod{\xi^2}$. Similarly, $Tf  = f + \frac{1}{2}\Lambda^{kj}\frac{\p f}{\p x^j}\xi_k \pmod{\xi^2}$. Fix local coordinates $\{x^k\}$ on $M$ and their copy $\{\tilde x^k\}$ on $\widetilde M$. Set $s^k(x,\nu, \xi) = S(x^k)$ and $t^k(x,\nu, \xi) = T(\tilde x^k)$. Then $s^k = x^k + \frac{1}{2}\Lambda^{kl}\xi_l \pmod{\xi^2}$ and $t^k = x^k + \frac{1}{2}\Lambda^{lk}\xi_l \pmod{\xi^2}$. Since both $S$ and $T$ are homomorphisms with respect to the point-wise products in $C^\infty(M)$ and $C^\infty(\T,Z)$, we observe that for $f\in C^\infty(M)$ the formal series in $\xi$, $Sf$ and $Tf$, can be written as  $Sf = f(s)$ and $Tf = f(t)$, where $f(s)$ and $f(t)$ are the compositions of the Taylor expansion of $f$ at $x=(x^1,\dots, x^n)$ (treated as a formal series) with the formal series $s= (s^1,\dots,s^n)$ and $t=(t^1,\dots, t^n)$, respectively.
The Poisson mapping $S\otimes T$ is induced by the formal mapping $s\times t:(x,\xi)\mapsto (s(x,\xi), t(x,\xi))$. Notice that $t^k - s^k = \frac{1}{2}\omega^{kl}\xi_l \pmod{\xi^2}$. Therefore (see the Appendix) if $M$ is symplectic then the mapping $s\times t$ is a formal symplectic diffeomorphism of $(\T, Z)$ onto $(M\times \widetilde M, M_{\rm diag})$ such that its restriction to the zero section $Z$ is the identity morphism onto the diagonal $M_{\rm diag}$ (both $Z$ and $M_{\rm diag}$ are copies of $M$). Since $M_{\rm diag}$ is a Lagrangian manifold in $M\times \widetilde M$,  it is well known that there exists a symplectic diffeomorphism of a tubular neighborhood of $M_{\rm diag}$  to a tubular neighborhood of $Z$ in $\T$ which is the identity morphism from $M_{\rm diag}$ to $Z$. Each local symplectic isomorphism of that sort induces a formal symplectic isomorphism between $(M\times \widetilde M, M_{\rm diag})$ and  $(\T, Z)$. The goal of this paper is to exhibit the particular formal symplectic diffeomorphism between $(M\times \widetilde M, M_{\rm diag})$ and  $(\T, Z)$ delivered by Fedosov's quantization.

{\it Remark.}  Notice that if $M$ is symplectic, $M\times \widetilde M$ is a (pair) symplectic groupoid. Thus both $(M\times \widetilde M, M_{\rm diag})$ and  $(\T, Z)$ can be thought of as isomorphic formal symplectic groupoids.

\section{Fedosov Star Products}

Fedosov's construction of star-products can be trivially generalized to the following setting (see \cite{BW}, \cite{DLS}, \cite{KSch},\cite{N}). Let $\Lambda^{jk}$ be a global tensor field on a symplectic manifold $(M,\omega)$ such that $\omega^{jk} = (1/2)(\Lambda^{jk} - \Lambda^{kj})$ is the nondegenerate Poisson tensor on $M$ inverse to the symplectic form $\omega$. Assume that $\nabla$ is an affine connection (possibly with torsion) which respects $\Lambda^{jk}$, i.e. $\nabla \Lambda^{jk} =0$. For each $x\in M$ we associate to the tangent space $T_xM$ a formal associative algebra $W_x$ whose elements are formal series
\[
                  a(\nu,y) = \sum_{r\geq 0, \alpha} \nu^r a_{r,\alpha} y^\alpha,
\]
where $\{y^k\}$ are linear coordinates on $T_xM$, $\alpha$ is a multi-index, and the standard multi-index notation is used.
The product in $W_x$ is given by the formula

\begin{equation}\label{E:prod}
   (a\circ b)(\nu, y)  = \exp \left(\frac{i\nu}{2}\Lambda^{jk}(x) \frac{\p^2}{\p y^j \p z^k}\right)a(\nu,y)b(\nu,z)|_{z=y}. 
\end{equation}
Taking the union of algebras $W_x$ we obtain a bundle $W$ of algebras. The fibre product (\ref{E:prod}) is extended to the sections of the bundle $W\otimes \Lambda$ of $W$-valued differential forms on $M$ by means of the usual exterior product of differential forms.
Introduce the gradings $\deg_\nu,\deg_a,\deg_s$ on the sections of $W\otimes \Lambda$ by defining it on the generators as follows:
\[
   \deg_\nu (\nu) =1, \quad \deg_s(y^k) = 1, \quad \deg_a (dx^k) = 1.
\]
All other gradings of the generators are set to zero. The product $\circ$ on the sections of $W\otimes \Lambda$ is bigraded with respect to the total grading ${\rm Deg} = 2\deg_\nu + \deg_s$ and the grading $\deg_a$. The connection $\nabla$ is extended to the sections of $W$ as follows:
\[
  \nabla a = \left(\frac{\p a}{\p x^j} - \Gamma^l_{jk} y^k \frac{\p a}{\p y^l}\right)dx^j,
\]
where $\Gamma^l_{jk}$ are the Christoffel symbols of $\nabla$.
Introduce the Fedosov operators $\delta$ and $\delta^{-1}$ on the sections of $W\otimes \Lambda$. Assume that a section $a\in \Gamma(W\otimes \Lambda)$ is homogeneous w.r.t. all the gradings introduced above with $\deg_sa=p$ and $\deg_aa = q$. Set
\begin{equation*}
   \delta(a) = dx^j \wedge \frac{\p a}{\p y^j} \text{\quad and 
\quad} \delta^{-1}a = 
\begin{cases}    
 \frac {1}{p+q}
 y^j i\left(\frac{\p}{\p x^j}\right) a 
   &\text{\rm if  $p+q > 0$,}
\\
 0, &\text{\rm if $p=q=0$.}
\end{cases}
\end{equation*}

Define the elements
\[
T:= \frac{1}{2}\, \omega_{s \alpha}T^\alpha_{kl} y^s dx^k 
\wedge 
dx^l
\mbox{\quad and \quad}
R := \frac{1}{4}\, \omega_{s\alpha} R^\alpha_{tkl}y^s y^t 
dx^k\wedge 
dx^l
\]
in $\Gamma(W\otimes \Lambda)$, where $T_{kl}^j = \Gamma_{kl}^j - \Gamma_{lk}^j $ is the torsion and
\[
 R^s_{tkl} := \frac{\p \Gamma^s_{lt}}{\p x^k} - \frac{\p 
\Gamma^s_{kt}}{\p x^l} + \Gamma^s_{k \alpha}\Gamma^\alpha_{lt} 
- \Gamma^s_{l \alpha}\Gamma^\alpha_{kt}
\]
is the curvature of the connection $\nabla$.

The following two theorems are minor modifications of the 
standard statements of Fedosov's theory adapted to the case of 
affine connections with torsion. We shall denote the Deg-homogeneous component of degree $k$ of an element $a\in \Gamma(W\otimes \Lambda)$ by $a^{(k)}$ and its $\deg_s$-homogeneous component of degree $l$ by $a_l$.
\begin{theorem}\label{T:fedcon} Let $\Omega = \sum_{k =1}^\infty \nu \Omega_k$ be a closed $\nu$-formal two-form on $M$. There exists a unique element $r \in \Gamma(W\otimes \Lambda)$ such that 
$r^{(0)}=r^{(1)}=0,\ \deg_a(r)=1,\ \delta^{-1}r = 0$, 
satisfying the equation 
\begin{equation}\label{E:delr}
\delta r = T + R + \nabla r - 
\frac{i}{\nu}\; r\circ r\ + \Omega.
\end{equation}
It can be calculated recursively with respect to the total 
degree Deg.  Then the Fedosov connection $D := -\delta + \nabla - 
\frac{i}{\nu} 
ad_\circ (r)$ is flat, i.e., $D^2=0$.
\end{theorem}
\noindent
The proof of the theorem is by induction, with the use of the 
identities
\[
\delta T = 0 \mbox{\quad and \quad} \delta R = \nabla T.
\]
which can be derived from the fact that the connection $\nabla$ respects the form $\omega$. 

The Fedosov connection $D$ is a $\deg_a$-graded derivation of 
the algebra $\Gamma(W\otimes \Lambda)$. Therefore Fedosov's algebra $\W_D := \ker D \cap \Gamma(W)$ is a subalgebra of $(\Gamma(W),\circ)$.

\begin{theorem}\label{T:fedquant}
The projection of $\W_D$  onto the 
part of $\deg_s$-degree zero, $w \mapsto w|_{y=0}$, is a bijection of Fedosov's algebra $\W_D$ onto $C^\infty(M)[[\nu]]$. The inverse 
mapping $\tau: C^\infty(M)[[\nu]] \to \W_D$ for a function $f 
\in C^\infty(M)$ can be calculated recursively w.r.t. the 
total degree Deg as follows:
\begin{gather*}
      \tau(f)^{(0)} = f,\\
      \tau(f)^{(k+1)} = \delta^{-1} \left(\nabla \tau(f)^{(k)} 
- \frac{i}{\nu} \sum_{l=0}^{k-1} 
ad_\circ\bigl(r^{(l+2)}\bigr)\bigl(\tau(f)^{(k-l)}\bigr)
  \right), k\geq 0.
\end{gather*}
The product $\ast$ on $C^\infty(M)[[\nu]]$ defined by the 
formula 
\[
f \ast g := (\tau(f) \circ \tau(g))|_{y=0}\ ,
\]
is a star-product on $M$.
\end{theorem}
\noindent Thus for $f \in C^\infty(M)[[\nu]]$ the element $\tau(f)\in \Gamma(W)$ is determined by the conditions $D\tau(f) = 0$ and $\tau(f)|_{y=0}=f$.
The following lemma and theorem can be proved as in \cite{BW} or \cite{N}.

\begin{lemma}\label{L:order}
For all $k \geq 1$ and $0 \leq l \leq \left[\frac{k-1}{2}\right]$ the mapping
\[
    C^\infty(M) \ni f \mapsto \tau(f)^{(k)}_{k-2l}
\]
is a differential operator of order $k-l$.
\end{lemma} 

Lemma \ref{L:order} implies

\begin{theorem}
   The Fedosov star-product $\ast$ is natural.
\end{theorem}

For different choices of the tensor $\Lambda^{jk}$ and connection $\nabla$ the star-product $\ast$ delivers  Fedosov quantizations \cite{F1}, deformation quantizations with separation of variables (or of the Wick type) \cite{BW}, \cite{CMP1}, \cite{N}, and the almost-K\"ahler deformation quantization \cite{KSch}.

\section{Formal geometric version of Fedosov's construction}

In this section we describe the formal geometric structure obtained from Fedosov's quantization by ``setting $\nu$ to zero".
This structure was discovered in \cite{EW} and later used in \cite{GL} as the first step of the BRST quantization procedure that yields Fedosov's quantization (conversion of the second class constraints to the first class constraints). 

For a section $w = w(x,\nu, y) \in \Gamma(W)$  denote by $w^{\vee} =  w|_{\nu = 0}$ the $\nu$-free part of $w$. Denote by $W^{\vee}$ the $\nu$-free part of the bundle $W$. Each fibre $W^{\vee}_x = \C[[y]]$ is an algebra with respect to the pointwise product. The smooth sections of $W^{\vee}$ can be canonically identified with the elements of $C^\infty(TM,Z)$ (abusing notations we denote by $Z$ also the zero section of $TM$). For two sections $a,b \in \Gamma(W)$ we have from (\ref{E:star}) and (\ref{E:prod}):
\begin{equation}\label{E:cl}
       (a \circ b)^{\vee} = a^{\vee} \cdot b^{\vee} \text{\quad and \quad} \left(\frac{1}{i\nu}[a,b]_\circ\right)^{\vee} = \{a^{\vee},b^{\vee}\}_{TM},
\end{equation} 
where $\{\cdot,\cdot\}_{TM}$ is the fibrewise Poisson bracket on $TM$ given on $a,b \in C^\infty(TM)$ by the formula
\begin{equation}\label{E:clpoiss}
      \{a,b\}_{TM} = \omega^{jk}(x)\frac{\p a}{\p y^j}\frac{\p b}{\p y^k}.
\end{equation} 

\begin{proposition}\label{P:poistau}
The mapping $C^\infty(M) \ni f\mapsto \tau(f)^{\vee}$ is a Poisson morphism from $C^\infty(M)$ to $C^\infty(TM,Z)$ endowed with the fibrewise Poisson bracket (\ref{E:clpoiss}).
\end{proposition}
\begin{proof}
Since the mapping $C^\infty(M)[[\nu]] \ni f\mapsto \tau(f)$ is $\nu$-linear,
$\tau(f)^{\vee}$ depends only on the $\nu$-free part of $f,\ f^{\vee} = f|_{\nu =0}$. For $f,g\in C^\infty(M)$ we have from (\ref{E:cl}) that
$\tau(f\ast g)^{\vee} = \left(\tau(f) \circ \tau(g)\right)^{\vee} =\tau(f)^{\vee}\cdot \tau(g)^{\vee}.$ On the other hand, $f \ast g = fg \pmod{\nu}$, whence $\tau(f \ast g)^{\vee} = \tau(f \cdot g)^{\vee}$. Thus $\tau(f \cdot g)^{\vee} = \tau(f)^{\vee}\cdot \tau(g)^{\vee}.$ Since $\frac{1}{i\nu}[f,g]_{\ast} = \{f,g\} \pmod{\nu}$ we get from (\ref{E:cl}) that 
\[
\tau(\{f,g\})^{\vee} = \tau\left(\frac{1}{i\nu}[f,g]_{\ast}\right)^{\vee} = \left(\frac{1}{i\nu}\left[\tau(f),\tau(g)\right]_{\circ}\right)^{\vee} = \{\tau(f)^{\vee},\tau(g)^{\vee}\}_{TM}.
\]
\end{proof}

In local coordinates on $M$ set $\vk^k = \tau(x^k)^{\vee}$. It follows from Proposition \ref{P:poistau} that the mapping $C^\infty(M) \ni f\mapsto \tau(f)^{\vee}$ is induced by the dual (formal) Poisson morphism $\vk: (x,y) \mapsto \vk(x,y)$ from $(TM,Z)$ to $M$, that is $\tau(f)^{\vee}(x,y) = f(\vk(x,y))$ (which is understood as above as a composition of formal series). 

The following theorem is the $\nu$-free version of Theorem ~\ref{T:fedcon}.  It provides an independent way of calculating $r^{\vee}$ and introduces a ``$\nu$-free  Fedosov connection" $D^{\vee}$ on $W^{\vee}$.  Denote by $\Xi_a$ the formal Hamiltonian vector field on $TM$ corresponding to a section $a$ of $W^{\vee}$  with respect to the fibrewise Poisson bracket (\ref{E:clpoiss}). 

\begin{theorem}\label{T:clfedcon}  The $\nu$-free  part $r^{\vee}$ of the element $r$ can be described as a unique element $\Gamma(W^{\vee}\otimes \Lambda)$ such that 
$r^{\vee}_0=r^{\vee}_1=0,\ \deg_a(r^{\vee})=1,\ \delta^{-1}r^{\vee} = 0$, 
satisfying the equation 
\begin{equation}\label{E:cldelr}
\delta r^{\vee} = T + R + \nabla r^{\vee} +
\frac{1}{2} \{r^{\vee},r^{\vee}\}_{TM}.
\end{equation}
It can be calculated recursively with respect to the symmetric 
degree $\deg_s$ as follows:
\begin{gather*}
                 r^{\vee}_2 = \delta^{-1} T,\\
   r^{\vee}_3 = \delta^{-1}\left(R + \nabla r^{\vee}_2 + 
   \frac{1}{2}\{r^{\vee}_2,r^{\vee}_2\}_{TM}\right) ,\\
   r^{\vee}_{k+3} = \delta^{-1}\left(\nabla r^{\vee}_{k+2} + 
\frac{1}{2} \sum_{l=0}^{k} \{r^{\vee}_{l+2}, r^{\vee}_{k-l+2}\}_{TM} 
\right), k\geq 1.
\end{gather*}
Then the connection $D^{\vee} := -\delta + \nabla +
\Xi_{r^{\vee}}$ on $W^{\vee}$ is flat, i.e., $\left(D^{\vee}\right)^2=0$.
\end{theorem}

Equation (\ref{E:cldelr}) is obtained from (\ref{E:delr}) with the use of (\ref{E:cl}) and an obvious identity $r \circ r = \frac{1}{2}[r,r]_\circ$. Notice that $r^{\vee}$ depends only on the $\nu$-free part $\omega$ of the Abelian curvature $\omega+\Omega$ of Fedosov's quantization.
For $w\in\Gamma(W)$ it is straightforward that $(Dw)^{\vee}=D^{\vee}w^{\vee}$. Since $D\tau(f)=0$ for any $f \in C^\infty(M)[[\nu]]$, it follows that $D^{\vee}\tau(f)^{\vee}=0$.
For $f \in C^\infty(M)$ the section $\tau(f)^{\vee}\in \Gamma(W^{\vee})$ can be found from the conditions  $D^{\vee}\tau(f)^{\vee}=0$ and $\tau(f)^{\vee}|_{y=0}=f$.
\begin{theorem}\label{T:clfedquant}
For  $f \in C^\infty(M)$ the element $\tau(f)^{\vee}$ can be calculated recursively w.r.t. the symmetric degree  $\deg_s$ as follows:
\begin{gather*}
      \tau(f)^{\vee}_0 = f,\\
      \tau(f)^{\vee}_{k+1} = \delta^{-1} \left(\nabla \tau(f)^{\vee}_k 
+ \sum_{l=0}^{k-1} 
\{r^{\vee}_{l+2},\tau(f)^{\vee}_{k-l}\}_{TM}
  \right), k\geq 0.
\end{gather*}
\end{theorem}

Since the mapping $\tau: f \mapsto \tau(f)$ from $C^\infty(M)[[\nu]]$
to $\Gamma(W)$ is given by a $\nu$-formal $\nu$-linear differential operator, the value of $\tau(f)$ at a given arbitrary point $x\in M$ is determined only by the $\nu$-formal $\infty$-jet of $f$ at $x$. In \cite{X} it was proved that the mapping $\tau$ is actually a bijection of the space of $\nu$-formal $\infty$-jets of functions at $x$ onto the fibre $W_x$. Below we will prove a ``formal geometric" version of this statement which is equivalent to the statement itself.
\begin{proposition}\label{P:class}
For an arbitrary point $x\in M$ the mapping $\tau^{\vee}_x: f \mapsto \tau(f)^{\vee}(x)$ from $C^\infty(M)$ to the fibre $W^{\vee}_x$ establishes a bijection between the space of $\infty$-jets of functions at $x$ and $W^{\vee}_x$.
\end{proposition}
\begin{proof}
The statement of the proposition is equivalent to the statement that the restriction of the dual (formal) mapping $\vk: (TM,Z) \to M$ to the fibre $W^{\vee}_x$ is a formal diffeomorphism between $(W^{\vee}_x,0)$ and $(M,x)$. The latter statement follows from the fact that $\vk^k = x^k + y^k \pmod{y^2}$ which can be obtained from Theorem \ref{T:clfedquant}.
\end{proof}

\section{Dequantization of a Fedosov star-product}

Introduce a bilinear pairing $\langle \cdot,\cdot \rangle$ on the bundle $W$ as follows: for $a,b \in \Gamma(W)$ set $\langle a,b\rangle = (a \circ b)|_{y = 0}$. Thus, in particular, $f\ast g = \langle \tau(f),\tau(g)\rangle$.

For a section $w \in \Gamma(W)$ introduce $\nu$-formal differential operators $L[w]$ and $R[w]$ on $M$ such that for $f \in C^\infty(M)[[\nu]]\ L[w]f = \langle w, \tau (f)\rangle$ and $R[w]f = \langle \tau (f), w\rangle$. Thus $L_f = L[\tau(f)]$ and $R_f = R[\tau(f)]$.

\begin{proposition}\label{P:natural}
 For a section $w \in \Gamma(W)$ the operators $L[w]$ and $R[w]$ are natural. Their $\sigma$-symbols $\sigma(L[w])$ and $\sigma(R[w])$ depend only on the $\nu$-free part of $w, \ \sigma(L[w]) = \sigma(L[w^{\vee}]),\ \sigma(R[w]) = \sigma(R[w^{\vee}])$. 
\end{proposition}
\begin{proof}
We will prove only the statement of the proposition related to the operator $L[w]$.
Fix $f\in C^\infty(M)$ and consider the contribution of the $\deg_\nu$- and $deg_s$-homogeneous components of $w$ and $\tau(f)$,
$w^{(j)}_{j-2k}$ and $\tau(f)^{(l)}_{l-2m}$, to $L[w](f)=\langle w, \tau(f)\rangle$. Notice that $\deg_\nu w^{(j)}_{j-2k} = k$ and $\deg_\nu \tau(f)^{(l)}_{l-2m} = m$. It follows from (\ref{E:prod}) that this contribution is nonzero only if $j-2k = l-2m$. Then it equals
\begin{equation}\label{E:contrib}
     \frac{\nu^{l-2m}}{(l-2m)!}\sum_{|J|=|K|=l-2m} \Lambda^{JK} \left(\frac{\p}{\p y}\right)_J w^{(j)}_{j-2k} \left(\frac{\p}{\p y}\right)_K \tau(f)^{(l)}_{l-2m},
\end{equation}
where $J$ and $K$ are multi-indices. The $\deg_\nu$-grading of  term (\ref{E:contrib}) is $l-2m + k + m = l-m+k$, while the order of the operator $C^\infty(M) \ni f \mapsto \tau(f)^{(l)}_{l-2m}$ is $l-m$ according to Lemma \ref{L:order}. Thus any component (\ref{E:contrib}) of the operator $L[w]$ is natural and has a nontrivial $\sigma$-symbol only if $k=0$, whence the proposition follows.
\end{proof}

\begin{proposition}\label{P:lmult}
  The mappings $w \mapsto \sigma(L[w])$ and $w \mapsto \sigma(R[w])$ from $\Gamma(W)$ to $C^\infty(\T,Z)$ are homomorphisms with respect to the point-wise products.
\end{proposition}
\begin{proof}
We need to show that $ \sigma(L[w_1w_2]) = \sigma(L[w_1]) \sigma(L[w_2])$ for any $w_1,w_2\in \Gamma(W)$. Fix an arbitrary  point $x\in M$. It follows from Proposition \ref{P:class} 
that there are functions $f,g\in C^\infty(M)$ such that $\tau(f)^{\vee}(x) = w^{\vee}_1(x),\ \tau(g)^{\vee}(x) = w^{\vee}_2(x)$. The distribution $L[w]_x:\phi \mapsto \langle w,\tau(\phi)\rangle(x)$ is natural and depends only on the value $w(x)$, therefore its $\sigma$-symbol depends only on $w^{\vee}(x)$. Thus $(Sf)(x) = \sigma(L_f)(x) = \sigma(L[\tau(f)])(x) = \sigma(L[w_1])(x)$ and $(Sg)(x) = \sigma(L[w_2])(x)$. We have $\tau(f\ast g)^{\vee}(x) = (\tau(f)\circ \tau(g))^{\vee}(x) = \tau(f)^{\vee}(x) \cdot \tau(g)^{\vee}(x) = w^{\vee}_1(x) \cdot w^{\vee}_2(x)$, therefore $S(f\ast g)(x) = \sigma(L_{f\ast g})(x) = \sigma(L[\tau(f\ast g)^{\vee}](x)  = \sigma(L[w^{\vee}_1\cdot w^{\vee}_2])(x) = \sigma(L[w_1w_2])(x)$. Since
$f\ast g = fg \pmod{\nu}$ we get that $S(f\ast g) = S(fg)=Sf \cdot Sg$. Similarly, since $w_1\circ w_2 = w_1 \cdot w_2 \pmod{\nu}$ we have $\sigma(L[w_1 \circ w_2]) =\sigma(L[w_1 \cdot w_2])$. Thus $\sigma(L[w_1])(x)\cdot\sigma(L[w_2])(x) = (Sf)(x)\cdot (Sg)(x) = S(fg)(x) =  \sigma(L[w_1 \cdot w_2])(x)$ for any $x\in M$. The proof of the statement  concerning $w \mapsto \sigma(R[w])$ is similar.
\end{proof}

Fix an arbitrary coordinate neighborhood $(U, \{x^k\})$ on $M$. As usual, $\{y^k\}$ are the corresponding fibre coordinates on $TU$.
Introduce a $\nu$-formal differential operator $Z_p$ on $U$  by the formula
\begin{equation}\label{E:Zp} 
   Z_pf = i\nu\frac{\p\tau(f)}{\p y^p}|_{y=0}.
\end{equation}
For $w = y^p$ we get from (\ref{E:prod}) and (\ref{E:Zp}) that 
\begin{equation}\label{E:ll}
L[y^p]f = \langle y^p, \tau(f) \rangle = \left( y^p \tau(f) + \frac{i\nu}{2}\Lambda^{pq}\frac{\p \tau(f)}{\p y^q}\right)|_{y=0} = \frac{1}{2} \Lambda^{pq}Z_q
\end{equation}
and 
\begin{equation}\label{E:rl}
R[y^p]f = \langle \tau(f), y^p \rangle = \left( y^p \tau(f) + \frac{i\nu}{2}\Lambda^{qp}\frac{\p \tau(f)}{\p y^q}\right)|_{y=0} = \frac{1}{2} \Lambda^{qp}Z_q.
\end{equation}
Since $\omega^{pq} = \frac{1}{2}\left(\Lambda^{pq} - \Lambda^{qp}\right)$ we get from  (\ref{E:ll}) and (\ref{E:rl}) that
$L[y^p] - R[y^p] = \omega^{pq} Z_q$ and therefore $Z_q = \omega_{qp}(L[y^p] - R[y^p])$, where $\left(\omega_{qp}\right)$ is inverse to $\left(\omega^{pq}\right)$. It follows from Proposition \ref{P:natural} that $Z_p$ is a natural operator. We denote by $\zeta_p = \zeta_p(x,\xi) \in C^\infty(U)[[\xi]]$  the $\sigma$-symbol of $Z_p,\ \zeta_p = \sigma(Z_p)$. Since for $f \in C^\infty(U)$ the expansion of $\tau(f)$ into Deg-homogeneous components is of the form
\[
\tau(f) = f(x) + \frac{\p f}{\p x^k} y^k + \dots,
\]
we obtain that $Z_p = i\nu \frac{\p}{\p x^p} \pmod {\nu^2}$, whence $\zeta_p = \xi_p \pmod{\xi^2}$. 

The following proposition is a direct consequence of Proposition \ref{P:lmult} and formulas (\ref{E:ll}) and (\ref{E:rl}).

\begin{proposition}\label{P:subst}
For $w \in C^\infty(U)[[\nu,y]]$ the $\sigma$-symbols of the operators $L[w]$ and  $R[w]$ are given by the formulas
\begin{equation}\label{E:subst}
        \sigma(L[w])= w^{\vee} \left(x, \frac{1}{2} \Lambda^{\cdot j} \zeta_j\right) \mbox{ {\rm and} }  \sigma(R[w])= w^{\vee} \left(x, \frac{1}{2} \Lambda^{j\cdot} \zeta_j\right).
\end{equation}
\end{proposition}

\begin{corollary}
For $f,g \in C^\infty(U)$
\[
Sf = \left(\tau(f)^{\vee}\right)\left(x, \frac{1}{2} \Lambda^{\cdot j} \zeta_j\right)  \mbox{ {\rm and} } Tf = \left(\tau(f)^{\vee}\right)\left(x, \frac{1}{2} \Lambda^{j\cdot} \zeta_j\right).
\]
\end{corollary}
\begin{proof} 
Since $L_f = L[\tau(f)]$, we obtain from Proposition \ref{P:subst} that $Sf = \sigma(L_f) = \sigma(L[\tau(f)]) = \left(\tau(f)^{\vee}\right)\left(x, \frac{1}{2} \Lambda^{\cdot j} \zeta_j\right)$. The formula for $T$ can be proved similarly.
\end{proof}

The element $\zeta_p= \zeta_p(x,\xi)$ is a formal series in $\xi$. Since $\zeta_p = \xi_p \pmod{\xi^2}$, one can express $\xi_p$ as a formal series in the variables $\zeta$.  It turns out that there is a simple formula for $\xi_p = \xi_p(x,\zeta)$. From that formula the formal series $\zeta_p= \zeta_p(x,\xi)$ can be recovered.

\begin{theorem}\label{T:xi} 
The following formula expresses $\xi_p = \xi_p(x,\zeta)$ in terms of the tensor $\Lambda^{jk}$ and the $\nu$-free  part of Fedosov's $\nu$-formal 1-form $r = r_p dx^p$:
\begin{equation}\label{E:thm}
      \xi_p = \zeta_p - r^{\vee}_p\left(x, \frac{1}{2}\Lambda^{\cdot j} \zeta_j\right) +  r^{\vee}_p\left(x, \frac{1}{2}\Lambda^{j\cdot} \zeta_j\right).
\end{equation}
\end{theorem}
\begin{proof} For $f \in C^\infty(U)$ we have that $D\tau(f)=0$ or
\begin{equation}\label{E:dtauf}
     -\delta \tau(f) + \nabla \tau(f) - \frac{i}{\nu} [r, \tau(f)]_\circ = 0.
\end{equation}
In components (\ref{E:dtauf}) can be rewritten as follows:
\begin{equation}\label{E:comp}
     -\frac{\p \tau(f)}{\p y^p}  +  \frac{\p \tau(f)}{\p x^p} - \Gamma_{pq}^t y^q \frac{\p \tau(f)}{\p y^t} - \frac{i}{\nu}\big( r_p \circ \tau(f) - \tau(f) \circ r_p\big) = 0.
\end{equation}
Now multiply (\ref{E:comp}) by $i\nu$ and set $y=0$. Using the fact that $\tau(f)|_{y=0} = f$ we get
\begin{equation}\label{E:oper}
     - Z_p f  +  i\nu \frac{\p f}{\p x^p}  + \langle r_p, \tau(f)\rangle - \langle \tau(f), r_p\rangle = 0.
\end{equation}
Equation (\ref{E:oper}) means that 
\begin{equation}\label{E:opzero}
- Z_p  +  i\nu \frac{\p}{\p x^p}  + L[r_p] - R[r_p] =0.
\end{equation}
All the summands on the left hand side of (\ref{E:opzero}) are natural operators. Passing to their $\sigma$-symbols with the use of Proposition \ref{P:subst} we obtain formula (\ref{E:thm}).
\end{proof}
The morphism $s \times t: (\T,Z) \to (M\times \widetilde M, M_{\rm diag})$ is defined locally as follows: $s^k = \vk^k\left(x, \frac{1}{2}\Lambda^{\cdot j}\zeta_j\right),\ t^l = \vk^l\left(x, \frac{1}{2}\Lambda^{j\cdot}\zeta_j\right)$, where $\vk^k = \left(\tau(x^k)\right)^{\vee}$. It follows from the results derived above from Fedosov's quantization that $s \times t$ is a formal symplectic isomorphism. This can be checked directly in the ``formal geometric setting". The symplectic structure on  $M\times\widetilde M$ corresponding to the Poisson structure is given in local coordinates $\{x^k, \tilde x^l\}$  by the formula $\omega_{jk}(x)dx^j\wedge dx^k - \omega_{jk}(\tilde x)d\tilde x^j\wedge d\tilde x^k$, where $\omega_{ij}$ is inverse to $\omega^{jk}$. The symplectic form on $\T$ corresponding to the Poisson bracket $\{\cdot,\cdot\}_\T$ is $2dx^p\wedge d\xi_p$. It remains to show that
\[
  \omega_{jk}(s)ds^j\wedge ds^k - \omega_{jk}(t)d t^j\wedge dt^k = 2dx^p\wedge d\xi_p,
\]
where $s^k = s^k(x,\zeta),\ t^l = t^l(x,\zeta)$ and $\xi_p = \xi_p(x,\zeta)$ is given by Theorem~\ref{T:xi}. The check is based on the formula 
\[
   \frac{\p \vk^p(x,y)}{\p y^k} \omega_{pq}(\vk(x,y)) \frac{\p \vk^q(x,y)}{\p y^l} = \omega_{kl}(x), 
\]
which can be derived from the fact that  the mapping $(x,y) \mapsto \vk(x,y)$ is a Poisson morphism from $(TM,Z)$ to $M$.

\section{Appendix}

In the Appendix we briefly introduce terminology related to formal neighborhoods and their morphisms.

Let $X$ be an $N$-dimensional real manifold and $Y$ its $n$-dimensional submanifold. Denote by $I$ the ideal in the algebra $C^\infty(X)$ of  functions vanishing on $Y$. Denote $I^\infty = \bigcap_{k = 1}^\infty I^k$. The algebra of jets $C^\infty(X)/I^\infty$ can be thought of as ``the algebra of smooth functions on the formal neighborhood of $Y$ in $X$". We use the notation $C^\infty(X,Y) := C^\infty(X)/I^\infty$. Locally this algebra is modelled as follows. Let $(U,\{x^k\})$ be a coordinate neighborhood in $X$ such that $U \cap Y$ is given by the equations $x^{n+1}=0,\dots, x^N=0$. Then $C^\infty(U, U\cap Y)$ is isomorphic to the algebra of formal series $C^\infty(U\cap Y)[[x^{n+1},\dots,x^N]]$. If $f\in C^\infty(U)$ is a representative of an element in $C^\infty(U,U\cap Y)$, the corresponding formal series in 
$C^\infty(U\cap Y)[[x^{n+1},\dots,x^N]]$ is the Taylor expansion of $f(x)$ at the points of $U\cap Y$ with respect to the (formal) parameters $x^{n+1},\dots,x^N$. 

If $E$ is a vector bundle with the zero section $Z$, the algebra $C^\infty(E,Z)$ can be naturally identified with the formal series of the form $\phi_0 + \phi_1+\dots$, where $\phi_r\in C^\infty(E)$ is a homogeneous polynomial of degree $r$ on each fibre. 

A local model for a formal neighborhood of the diagonal in the Cartesian square of an $n$-dimensional manifold $M, \ C^\infty(M\times M, M_{\rm diag})$, can be obtained by taking two copies of coordinates $\{s^k\}$ and $\{t^k\}$ on a neighborhood $U \subset M$ so that the diagonal $(U\times U) \cap M_{\rm diag}$ is given by the equations $t^k - s^k =0$. Then we can use the coordinates $\{s^1, \dots, s^n, t^1 - s^1, \dots,t^n - s^n\}$ to model locally $C^\infty(M\times M,M_{\rm diag})$ as above.

A morphism $\Psi: C^\infty(X,P) \to C^\infty(Y,Q)$ can be described locally as follows. Let $U\subset X$ and $V \subset Y$ be two coordinate neighborhoods with the coordinates $(x,\eta)$ and $(y,\theta)$ respectively, where $U \cap P$ is defined by the equations $\eta=0$ and $V\cap Q$ is defined by $\theta =0$. Suppose that $x^k = x^k(y,\theta)$ and $\eta^l = \eta^l (y,\theta)$ are formal series in $\theta$ with coefficients in $C^\infty(V\cap Q)$ such that $\eta = 0\pmod{\theta}$. These data define a morphism $\Psi: C^\infty(U\cap P)[[\eta]] \to C^\infty(V\cap Q)[[\theta]]$ via the substitution $\Psi: f(x,\eta) \mapsto f(x(y,\theta), \eta(y,\theta))$ and subsequent Taylor expansion with respect to the formal parameters $\theta$. We can think of $\Psi$ as of a mapping induced by the mapping $\psi : (y,\theta) \mapsto (x(y,\theta),\eta(y,\theta))$. Let $x^k(y,\theta) = \phi^k (y) \pmod{\theta}$ and $\eta^m (y,\theta) = \chi^m_n (y)\theta^n \pmod{\theta^2}$. The mapping $\psi$ is invertible if $\phi:V\cap Q \to U\cap P$ is a diffeomorphism and the matrix $(\chi^m_n)$ is invertible.

\end{document}